\documentclass[12pt]{article}

\usepackage{amsmath}
\usepackage{amssymb}

\numberwithin{equation}{section}
\newcommand{\1}{\mathbf{1}}
\newcommand{\Z}{{\mathbb Z}}
\newcommand{\C}{{\mathbb C}}
\newcommand{\g}{{\mathfrak g}}

\newcommand{\VL}{V_{L}}

\newcommand{\widet}{\widetilde}

\newcommand{\al}{\alpha}
\newcommand{\be}{\beta}
\newcommand{\dl}{\delta}
\newcommand{\gm}{\gamma}
\newcommand{\om}{\omega}
\newcommand{\sg}{\sigma}

\newcommand{\prf}{\noindent {\bfseries Proof} \  }
\newcommand{\qed}{\mbox{ $\square$}}
\newcommand{\la}{\langle}
\newcommand{\op}{\oplus}
\newcommand{\ot}{\otimes}
\newcommand{\ra}{\rangle}

\newcommand{\Aut}{\mbox{Aut}}

\newtheorem{thm}{Theorem}[section]

\newtheorem{lem}[thm]{Lemma}

\begin{document}

\begin{center}
{\Large {\bfseries Decomposition of 
the vertex operator algebra $V_{\sqrt{2}A_{3}}$}}

\vspace{10mm}

Chongying Dong\footnote{Supported by NSF grant 
DMS-9700923 and a research grant from the Committee on Research, UC Santa Cruz.}\\
Department of Mathematics, University of California\\
Santa Cruz, CA 95064\\
\vspace{2mm}
Ching Hung Lam\\
Institute of Mathematics, University of Tsukuba\\
Tsukuba 305-8571, Japan\\
\vspace{2mm}
Hiromichi Yamada\\
Department of Mathematics, Hitotsubashi University\\
Kunitachi, Tokyo 186-8601, Japan

\end{center}
\vspace{10mm}


\section{Introduction}

A conformal vector with central charge $c$ in a vertex operator algebra 
is an element of weight two whose component operators
satisfy the Virasoro algebra relation with central
charge $c.$ Then the vertex operator subalgebra 
generated by the vector is isomorphic to a highest 
weight module for the Virasoro algebra with central charge $c$
and highest weight $0$ (cf. \cite{M}).

Let $V_{\sqrt{2}A_{l}}$ be the vertex operator algebra 
associated with a lattice $\sqrt{2}A_{l}$, where 
$\sqrt{2}A_{l}$ denotes $\sqrt{2}$ times an ordinary 
root lattice of type $A_{l}$. Motivated by the problem 
of looking for maximal associative algebras of the Griess
algebra [G], a class of conformal vectors in 
$V_{\sqrt{2}A_{l}}$ were studied and constructed in \cite{DLMN}. 
It was shown in \cite{DLMN} that the Virasoro element $\om$ 
of  $V_{\sqrt{2}A_{l}}$ is decomposed into a sum of $l+1$ mutually 
orthogonal conformal vectors $\om^{i}$; $1 \le i \le l+1$ 
with central charge $c_{i} = 1 - 6/(i+2)(i+3)$ for 
$1 \le i \le l$ and $c_{l+1} = 2l/(l+3)$.  The vertex
operator subalgebra generated by conformal vector $\om^i$
is exactly the irreducible highest weight module
$L(c_i,0)$ for the Virasoro algebra. The vertex operator
subalgebra $T=T_l$ generated by these conformal 
vectors is isomorphic to a tensor 
product $\ot_{i=1}^{l+1} L(c_{i},\, 0)$ of the Virasoro 
vertex operator algebras $L(c_{i},\, 0)$ and 
$V_{\sqrt{2}A_{l}}$ is a direct sum of irreducible 
$T$-submodules.

In this paper we determine the decomposition
of $V_{\sqrt{2}A_{3}}$ into the direct sum of irreducible
$T$-modules completely. 
The direct summands have been determined \cite{KMY} in the case 
$l=2.$ For general $l$ only the direct summands with minimal weights 
at most two are known \cite{Y}.

The main idea for the decomposition in this paper is to embed
$V_{\sqrt{2}A_{3}}$ into the vertex operator algebra 
$V_{(\sqrt{2}A_{1})^{\oplus 3}}$ by considering  
the lattice $\sqrt{2}A_{3}$ as a sublattice of 
$(\sqrt{2}A_{1})^{\oplus 3}.$ It turns out that
$V_{\sqrt{2}A_{3}}$ is isomorphic to the vertex
operator subalgebra $V_{(\sqrt{2}A_{1})^{\oplus 3}}^+$ 
which is the fixed points of the involution of
$V_{(\sqrt{2}A_{1})^{\oplus 3}}$ induced from the $-1$
isometry of $(\sqrt{2}A_{1})^{\oplus 3}.$ Moreover,  
$V_{(\sqrt{2}A_{1})^{\oplus 3}}^+$ has a subalgebra isomorphic to
$V_{\sqrt{2}A_{2}}^+\otimes V_F^+$ where $F$ is a 
rank one lattice spanned by an element of square length 6 and
$V_{(\sqrt{2}A_{1})^{\oplus 3}}^+$ is a direct sum of 3 irreducible
modules for $V_{\sqrt{2}A_{2}}^+\otimes V_F^+.$
Then using \cite{DG} on the decomposition
of lattice type vertex operator algebra of rank 1 
into the direct sum of irreducible modules for the Virasoro algebra 
and results in \cite{KMY}, we can determine all the irreducible
$T$-modules in  $V_{\sqrt{2}A_{3}}.$ 
We should also mention that the sum of certain irreducible
modules for $T$ inside $V_{\sqrt{2}A_{3}}$  forms a rational
vertex operator algebra from our picture. 

The authors would like to thank Masaaki Kitazume, Geoffrey 
Mason, and Masahiko Miyamoto for stimulating conversations.

\section{Some automorphisms of $V_{\Z\al}$}

Our notation for the vertex operator algebra $\VL = 
M(1) \ot \C[L]$ associated with a positive definite even 
lattice $L$ is standard \cite{FLM}. In particular,
${\mathfrak h}=\C\otimes_{\Z} L$  is an abelian Lie algebra
and extend the the bilinear form to ${\mathfrak h}$ by 
$\C$-linearity,
$\hat {\mathfrak h}={\mathfrak h}\otimes \C[t,t^{-1}]\oplus \C K$
is the corresponding affine algebra, 
$M(1)=\C[\al(n)|\al\in {\mathfrak h}, n<0],$ 
where $\al(n)=\al\otimes t^n,$ is the unique irreducible
   ${\mathfrak h}$-module such that $\alpha(n)1=0$ 
for all $\alpha\in {\mathfrak h}$ and $n$ positive,
and $K=1.$ The element in 
the group algebra $\C[L]$ of the additive group $L$ 
corresponding to $\beta \in L$ will be denoted by 
$e^{\beta}$. Note that the central extension $\hat L$ of $L$ 
by the cyclic group
of order 2 is split if the square length of any element
in $L$ is a multiple of 4 (cf. \cite{FLM}). For example,
$\sqrt{2}A_l$ is a such lattice. The vacuum vector $\1$  of $V_L$
is $1\otimes e^0$ and the Virasoro element $\omega$
is $\frac{1}{2}\sum_{i=1}^d\beta_i(-1)^2$ where $\{\beta_1,...,
\beta_d\}$ is an orthonormal basis of  ${\mathfrak h}.$

We need to know explicit expressions of the vertex operators $Y(u,z)$
for $u=h(-1)$ or $u=e^{\beta}$ for $h\in  {\mathfrak h}$ and $\beta\in L$
in the next section to do certain calculations. We assume that 
the square length of any element in $L$ is a multiple of 4. 
The operator $Y(h(-1),z)$ 
is defined as
\begin{equation}\label{ea1}
Y(h(-1),z)=\sum_{n\in\Z}h(n)z^{-n-1}=\sum_{n\in\Z}h(-1)_nz^{-n-1}
\end{equation}
where $h(n)$ acts on $M(1)$ if $n\ne 0$ and $h(0)$ acts on $\C[L]$
so that $h(0)e^{\gamma}=\la h,\gamma\ra e^{\gamma}$ for
$\gamma\in L.$ In order to define $Y(e^{\beta},z)$ we need
to define operators $e_{\beta}$ and $z^{\beta}$ acting on $V_L$
such that $e_{\beta} (u\otimes e^{\gamma})=u\otimes e^{\beta+\gamma}$
and $z^{\beta}(u\otimes e^{\gamma})=z^{\la \beta,\gamma\ra}u\otimes e^{\gamma}$
for $u\in M(1)$ and $\gamma\in L.$ Then 
\begin{equation}\label{ea2}
Y(e^{\beta},z)=\sum_{n\in \Z}e^{\beta}_nz^{-n-1}=
e^{\sum_{n<0}\frac{\beta(n)}{-n}z^{-n}}
e^{\sum_{n>0}\frac{\beta(n)}{-n}z^{-n}}e_{\beta}z^{\beta}.
\end{equation}
We refer the reader to \cite{FLM} for the definition of vertex
operators $Y(u,z)$ for general $u\in V_L.$ 

Let $L^{\circ}=\{\alpha\in {\mathfrak h}| \la \alpha, L\ra\subset \Z\}$
be the dual lattice of $L.$ Then $L^{\circ}/L$ is a finite group.
For each $\lambda\in L^{\circ}$ the corresponding untwisted Fock space
$V_{L+\lambda}=M(1)\otimes \C[L+\lambda]$ is an irreducible  module
for $V_L$ \cite{FLM}. Let $L^{\circ}=\cup_{i\in L^{\circ}/L}(L+\lambda_i)$
be a coset decomposition. Then $V_{L+\lambda_i}$ are all
inequivalent irreducible $V_L$-modules \cite{D}. 

Let $V_{\Z\al}$ be the vertex operator algebra 
associated with a rank one lattice $\Z\al$ such 
that $\la \al, \al \ra = 2.$ The homogeneous 
subspace $\g = (V_{\Z\al})_{(1)}$ of $V_{\Z\al}$ 
of weight one possesses a Lie algebra structure 
given by $[u, v] = u_{0}v$ with a symmetric 
invariant form $\la \cdot,\cdot\ra$ such that  
$\la u,v\ra \1 = u_{1}v$ (\cite[Section 8.9]{FLM}). 
We have
\begin{align*}
[e^{\al}, e^{-\al}] &= \al(-1), & [\al(-1),e^{\pm\al}] &=
\pm 2e^{\pm\al},\\
\la e^{\al}, e^{-\al}\ra &= 1, & \la \al(-1), \al(-1) \ra 
& = 2, 
\end{align*}
and $\la u, v \ra = 0$ for the other pairs $u, v$ in 
$\{ \al(-1), e^{\al}, e^{-\al}\}$.  In particular, 
$\{ \al(-1), e^{\al}, e^{-\al}\}$ is a standard basis of 
$\g \cong sl_{2}(\C)$.

Now consider three automorphisms (cf. \cite{FLM})
$\theta_{1}$, $\theta_{2}$, $\sg$ of $(\g, \la \cdot,\cdot\ra )$ 
of order two such that 
\begin{align*}
\theta_{1} & : \al(-1) \longmapsto \al(-1), & 
e^{\al}  & \longmapsto -e^{\al}, & 
e^{-\al} & \longmapsto -e^{-\al},\\
\theta_{2} & : \al(-1) \longmapsto -\al(-1), & 
e^{\al}  & \longmapsto e^{-\al}, & 
e^{-\al} & \longmapsto e^{\al},\\
\sg & : \al(-1) \longmapsto e^{\al} + e^{-\al}, & 
e^{\al} + e^{-\al} & \longmapsto
 \al(-1), & e^{\al} - e^{-\al} & \longmapsto -(e^{\al} - e^{-\al}). 
\end{align*}
 
Clearly $\sg \theta_{1} \sg = \theta_{2}$.  These automorphisms 
of $(\g, \la \cdot,\cdot\ra )$ can be uniquely extended to 
automorphisms of the vertex operator algebra $V_{\Z\al}$. 
In order to see this we recall that the Verma module $V(1,0)$ for
the affine algebra $A_1^{(1)}=sl_2(\C)\otimes \C[t,t^{-1}]\oplus \C {\bf c}$
is the quotient
of $U(A_1^{(1)})$ modulo the left ideal generated
by $x\otimes t^n,$ ${\bf c}-1$ for $x\in sl_2(\C)$ and
$n\geq 0.$ Note that the  automorphism
group $\Aut(sl_2(\C))$ of the Lie algebra $sl_2(\C)$  acts on $A_1^{(1)}$ by
acting on the first tensor factor of $sl_2(\C)\otimes \C[t,t^{-1}]$
and trivially on ${\bf c}.$ As a result $\Aut(sl_2(\C))$ acts
on $U(A_1^{(1)})$ as algebra automorphisms. Clearly
this induces an action of $\Aut(sl_2(\C))$  on $V(1,0).$ 
Note that $V_{\Z\alpha}$
is the irreducible quotient of $V(1,0)$ modulo the maximal 
submodule for $A_1^{(1)}.$ It is easy to see from this  construction
that $\Aut(sl_2(\C))$ acts on $V_{\Z\alpha}.$ In fact, the subgroup
of $\Aut(sl_2(\C))$ consisting of those preserving the invariant bilinear
form can be regarded as a subgroup of the automorphisms of the vertex
operator algebra $V_{\Z\alpha}.$ 
(This observation works for any finite dimensional
semisimple Lie algebra in the position of $sl_2(\C).$) 
Since $\theta_1, \theta_2$ and $\sigma$ preserve the bilinear form
on $sl_2(\C)$ they act on  $V_{\Z\alpha}$ as vertex operator algebra 
automorphisms.

We denote the corresponding automorphisms of $V_{\Z\al}$ by the same symbols 
$\theta_{1}$, $\theta_{2}$, and $\sg$.  Then on $V_{\Z\al}$, we have 
$\theta_{1}(u \ot e^{\beta}) = (-1)^{\la \al, 
\beta \ra/2} u \ot e^{\beta}$ for $u\in M(1)$ and $\beta \in \Z\al$ and 
$\theta_{2}$ is the automorphism induced from the isometry 
$\beta \longmapsto -\beta$ of $\Z\al$ \cite{FLM}. We still have 
$\sg \theta_{1} \sg 
= \theta_{2}$. We also have
\[
\sg(\al(-1)^{2}) = \al(-1)^{2} \quad \mbox{and} \quad 
\sg(e^{\pm \al}) = \frac{1}{2}(\al(-1) \mp (e^{\al} - e^{-\al})).
\]
We should mention that $\sg(\al(-1)^{2}) = \al(-1)^{2}$ is not obvious.
Using the definitions of $Y(e^{\pm \alpha},z)$ and $\sg$ we see that
\begin{align*}
\sg(\al(-1)^{2}) &= \sg(\al(-1))_{-1}\sg(\al(-1))\\
&= (e^{\alpha}+e^{-\alpha})_{-1}(e^{\alpha}+e^{-\alpha})\\
&= \al(-1)^{2}.
\end{align*}

\section{Decomposition of $V_{\sqrt{2}A_{3}}$}

Let $L$ be a lattice with basis 
$\{\al_{1},\al_{2},\al_{3}\}$ such that $\la \al_{i},\al_{j}\ra 
= 2\dl_{ij}$.  Then $L=A_1\oplus A_1\oplus A_1$  where
$A_1$ is the root lattice of $sl_2(\C).$ Set
\[
\beta_{1}=(\al_{1}+\al_{2})/\sqrt{2}, \qquad 
\beta_{2}=(-\al_{2}+\al_{3})/\sqrt{2}, \qquad 
\beta_{3}=(-\al_{1}+\al_{2})/\sqrt{2}.
\]
Then $\{\beta_{1}, \beta_{2}, \beta_{3}\}$ forms the set of 
simple roots of type $A_{3}$.  Set $\gm = -\al_{1}+\al_{2}+\al_{3}$. 
We consider two sublattices of $L:$ 
\[
N = \sum_{i,j=1}^{3} \Z (\al_{i} \pm \al_{j}), \qquad 
D = E \op F,
\]
where $E = \Z (\al_{1}+\al_{2}) + \Z (-\al_{2}+\al_{3}) 
= \Z\sqrt{2}\beta_{1} + \Z\sqrt{2}\beta_{2}$ and $F = \Z\gm$. 

\begin{lem}\label{l3.1} $(1)$ We have that 
$N = \{ \beta \in L\,|\, \la \al_{1}+\al_{2}+\al_{3},\, \beta \ra 
\equiv 0 \pmod 4 \},$ $N$ is isometric to $\sqrt{2}A_{3},$ 
 and $[L:N]=2$.

$(2)$ $[L:D]=3$ and 
$L = D \cup (D+\al_{2}) \cup (D-\al_{2})$

$(3)$ $\al_{2} = \sqrt{2}(\beta_{1}-\beta_{2})/3 + \gm/3$ and 
each element of the coset $D+\al_{2}$ can be 
uniquely written as an 
orthogonal sum of an element in  
$E+\sqrt{2}(\beta_{1}-\beta_{2})/3$ and an element in  
$F+\gm/3$.

$(4)$ $E$ is isometric to $\sqrt{2} A_2.$ 
\end{lem}

\prf (1) The first assertion can be verified easily.  Since $N = 
\Z\sqrt{2}\beta_{1} + \Z\sqrt{2}\beta_{2} + 
\Z\sqrt{2}\beta_{3}$, the second assertion holds.  We have 
$\al_{i} \not\in N$ and $L = N \cup (N+\al_{i})$ for any $i$. 
In particular $[L:N]=2$. 

(2)--(4) are obvious.\quad \qed

\medskip

Since $E$ and $F$ are even lattices $V_E$ and $V_F$ are vertex
operator algebras which can be regarded as vertex operator 
subalgebras of $V_{\sqrt{2}A_1^{\oplus 3}}$ (with different
Virasoro algebras). 
Note that $\la E+\sqrt{2}(\beta_{1}-\beta_{2})/3, E\ra\subset \Z$ 
and $\la F+\gamma/3,F\ra \subset \Z.$ Thus $V_{E+\sqrt{2}(\beta_{1}-\beta_{2})/3}$ is an irreducible $V_E$-module and $V_{F+\gamma/3}$ is
an irreducible $V_F$-module (cf. \cite{FLM}).

The lattice $L$ is a direct sum of $\Z\al_{i}$; $i=1$, $2$, $3$ 
and thus the vertex operator algebra $V_{L}$ associated with $L$ is 
a tensor product $V_{L} = V_{\Z\al_{1}} \ot 
V_{\Z\al_{2}} \ot V_{\Z\al_{3}}$ (see \cite{FHL} for the definition
of tensor product vertex operator algebra).
Define three automorphisms 
of $V_{L}$ of order two by 
\[
\psi_{1} = \theta_{1} \ot \theta_{1} \ot \theta_{1}, \qquad 
\psi_{2} = \theta_{2} \ot \theta_{2} \ot \theta_{2}, \qquad 
\tau = \sg \ot \sg \ot \sg,
\]
where $\theta_{1}$, $\theta_{2}$, and $\sg$ are the automorphisms 
of $V_{\Z\al_{i}}$ described in Section 2. Then 
\[
\psi_{1}(u \ot e^{\beta}) 
= (-1)^{\la \al_{1}+\al_{2}+\al_{3}, \,\beta \ra /2}u \ot e^{\beta}
\]
for $u \in M(1)$ and $\beta \in L$, $\psi_{2}$ is the automorphism 
induced from the isometry $\beta \longmapsto -\beta$ of $L$, and 
$\tau \psi_{1} \tau = \psi_{2}$.

For any $\psi_{2}$-invariant subspace $U$ of $V_{L}$ we shall 
write $U^{\pm} = \{ v \in V_{L}\,|\, \psi_{2}(v) = \pm v\}$. 

\begin{lem}\label{l3.2} 
We have $\tau (V_{N}) = V_{L}^{+}$.

\end{lem}

\prf From the action of $\psi_{1}$ on $V_{L}$ and Lemma \ref{l3.1} (1) 
it follows that $V_{N} = \{ v \in V_{L}\,|\, \psi_{1} (v) = 
v \}$.  Since $\psi_{2} \tau = \tau\psi_{1}$, the 
assertion holds.
\quad \qed 

\medskip
\begin{lem}\label{l3.3}
$(1)$ \  $V_{L}^{+} \cong V_{D}^{+} \op V_{D+\al_{2}}$ as 
$V_{D}^{+}$-modules.\\
$(2)$ \  $V_{D}^{+} = (V_{E}^{+} \ot V_{F}^{+}) \op 
(V_{E}^{-} \ot V_{F}^{-})$.\\
$(3)$ \  $V_{D+\al_{2}} = V_{E+\sqrt{2}(\beta_{1}-\beta_{2})/3} 
\ot V_{F+\gm/3}$.
\end{lem}

\prf  Lemma \ref{l3.1} (2) implies $V_{L} = V_{D} \op V_{D+\al_{2}} 
\op V_{D-\al_{2}}$.  Since $\psi_{2}$ leaves $V_{D}$ invariant 
and interchanges $V_{D+\al_{2}}$ and $V_{D-\al_{2}}$, 
we have
\begin{align*}
V_{L}^{+} &= V_{D}^{+} \op (V_{D+\al_{2}} \op V_{D-\al_{2}})^{+}\\
& \cong V_{D}^{+} \op V_{D+\al_{2}}
\end{align*}
as $V_{D}^{+}$-modules.  Now $D = E \op F$ and both of $E$ and 
$F$ are $\psi_{2}$-invariant.  Hence (2) holds.  (3) follows 
from Lemma \ref{l3.1} (3). \quad \qed

\medskip
Similarly we have the following result which is not used in this
paper to study the decomposition of $V_N$ into the sum of irreducible
$T$-modules.

\begin{lem}\label{l3.4}
$(1)$ \   $V_{L}^{-} \cong V_{D}^{-} \op V_{D+\al_{2}}$ as 
$V_{D}^{+}$-modules.\\
$(2)$ \  $V_{D}^{-} = (V_{E}^{+} \ot V_{F}^{-}) \op 
(V_{E}^{-} \ot V_{F}^{+})$.
\end{lem}

Let us recall the conformal vectors introduced in [DLMM]. 
For any positive integer $l>0$   let $\Phi_l$ be the 
root system of type $A_l$ generated by simple
roots $\beta_1,...,\beta_l$ and $\Phi_l^+$ be the set of
positive roots. 
We assume that the square length of
a root is 2. We also assume that $\Phi_i$ is a sub-system of
$\Phi_{l}$ with simple roots $\be_1,...\beta_i$ for $i\leq l.$ 
Let $\sqrt{2}A_i$ be the positive definite even lattice spanned
by $\sqrt{2} \Phi_i.$ Then $\sqrt{2}A_i$ is a sublattice 
of $\sqrt{2}A_l$  and $V_{\sqrt{2}A_i}$
is a vertex operator subalgebra of $V_{\sqrt{2}A_l}.$ 

For $\beta \in \Phi_l$ set 
\[
w_{\be}^{\pm} = \be(-1)^{2} \pm 2(e^{\sqrt{2}\be} + e^{-\sqrt{2}\be})
\]
and
\[
s^{i} = \frac{1}{2(i+3)}\sum_{\be \in \Phi_{i}^{+}} w_{\be}^{-},\qquad 
\om = \frac{1}{2(l+1)}\sum_{\be \in \Phi_{l}^{+}} \be(-1)^{2}.
\]
Then the conformal vectors of $V_{\sqrt{2}A_{l}}$ 
defined in \cite{DLMN} are
\begin{equation}\label{ea3}
\om^{1}=s^{1}, \qquad \om^{i+1}=s^{i+1}-s^{i}, i=1,...,l-1,  
\qquad \om^{l+1}=\om-s^{l}.
\end{equation}
Note that $\om$ is the Virasoro element.

By Lemma \ref{l3.2}, $V_{N}$ is isomorphic to $\VL^{+}$ as a 
vertex operator algebra.  We next calculate the images 
of the conformal vectors $\om^{1}$, $\om^{2}$, $\om^{3}$, 
and $\om^{4}$ in $V_{N}$ under  
the automorphism $\tau$.  Note that $\Phi_{2}^{+} = 
\{\beta_{1}, \beta_{2}, \beta_{1}+\beta_{2}\}$ and $\Phi_{3}^{+} 
= \Phi_{2}^{+} \cup 
\{\beta_{3}, \beta_{2} + \beta_{3}, \beta_{1}+\beta_{2}+\beta_{3}\}$. 

Let 
$L(c,h)$ be the irreducible highest weight module for the Virasoro
algebra with central charge $c$ and highest weight $h.$ Then 
$L(c,0)$ is a vertex operator algebra if $c\ne 0$ (cf. \cite{FZ}). 
The conformal vectors $\om^{1}$, $\om^{2}$, $\om^{3}$, 
and $\om^{4}$ 
are mutually orthogonal and their 
central charges are $1/2$, $7/10$, $4/5$, and $1$. 
Hence the 
subalgebra $T$ of $V_{N}$ generated by these conformal 
vectors is isomorphic to 
\[
  \begin{array}{l}
L(\frac{1}{2},\,0) \ot L(\frac{7}{10},\,0) \ot 
L(\frac{4}{5},\,0) \ot L(1,\, 0).
  \end{array}
\]
Moreover, $V_N$ is a completely reducible $T$-module. 
The main purpose in this paper is to determine the irreducible
$T$-modules in $V_N$ or equivalently in $V_L^+.$

In order to achieve this we need to determine the images of 
$s^i$ under $\tau.$ 
\begin{lem}\label{al1} We have
\begin{align*}
\tau(s^{1}) &= \frac{1}{8}w_{\be_{3}}^{+},\\
\tau(s^{2}) &= \frac{1}{10}(w_{\be_{3}}^{+} + w_{\be_{2}}^{-} 
+ w_{\be_{2}+\be_{3}}^{+}),\\
\tau(s^{3}) &= \frac{1}{12}(w_{\be_{3}}^{+} + w_{\be_{2}}^{-} 
+ w_{\be_{2}+\be_{3}}^{+} + w_{\be_{3}}^{-} + 
w_{\be_{2}+\be_{3}}^{-} + w_{\be_{2}}^{+})\\
&= \frac{1}{6}(\be_{2}(-1)^{2} + \be_{3}(-1)^{2} 
+ (\be_{2}+\be_{3})(-1)^{2}),\\
\tau(\om) &= \om.
\end{align*}
\end{lem}

\prf  The proof is a straightforward computation. We compute $\tau(s^1)$
here and leave the others to the reader. Note that
\begin{align*}
s^1 &=\frac{1}{8} \left(\be_1(-1)^2-2(e^{\sqrt{2}\be_1} + 
e^{-\sqrt{2}\be_1})\right)\\
 &=\frac{1}{16}(\alpha_1(-1)+\alpha_2(-1))^2 - 
\frac{1}{4}(e^{\alpha_1}e^{\alpha_2} + e^{-\alpha_1}e^{-\alpha_2}),
\end{align*}
where $e^{\pm\alpha_1}e^{\pm\alpha_2}$ is understood as tensor product
of  $e^{\pm\alpha_1}$ with $e^{\pm\alpha_2}$ in 
$V_{\Z\alpha_1}\otimes V_{\Z\alpha_2}.$ Recall the definitions
of vertex operators $Y(h(-1),z)$ and $Y(e^{\beta},z)$ from 
(\ref{ea1}) and (\ref{ea2}). Then from the definition
of $\tau$ we obtain
\begin{align*}
\tau(s^1) &=\frac{1}{16}(e^{\alpha_1}+e^{-\alpha_1}+
e^{\alpha_2}+e^{-\alpha_2})_{-1}
(e^{\alpha_1}+e^{-\alpha_1}+e^{\alpha_2}+e^{-\alpha_2})\\
&\qquad -\frac{1}{16}(\alpha_1(-1)-(e^{\alpha_1}-e^{-\alpha_1}))
(\alpha_2(-1)-(e^{\alpha_2}-e^{-\alpha_2}))\\
&\qquad -\frac{1}{16}(\alpha_1(-1)+(e^{\alpha_1}-e^{-\alpha_1}))
(\alpha_2(-1)+(e^{\alpha_2}-e^{-\alpha_2}))\\
&=\frac{1}{16}(\alpha_1(-1)^2+\alpha_2(-1)^2)+
\frac{1}{8}(e^{\alpha_1+\alpha_2}+e^{\alpha_1-\alpha_2}+e^{-\alpha_1+\alpha_2}
+e^{-\alpha_1-\alpha_2})\\
&\qquad -\frac{1}{8}\alpha_1(-1)\alpha_2(-1)-
\frac{1}{8}(e^{\alpha_1+\alpha_2}-e^{\alpha_1-\alpha_2}-e^{-\alpha_1+\alpha_2}
+e^{-\alpha_1-\alpha_2})\\
&=\frac{1}{8} (\be_3(-1)^2+2(e^{\sqrt{2}\be_3} + e^{-\sqrt{2}\be_3}))\\
&=\frac{1}{8}w_{\beta_3}^+.
\end{align*}
\qed

Clearly, $\tau(s^1), \tau(s^2-s^1), \tau(s^3-s^2)$ are not the conformal
vectors associated to the lattice $\sqrt{2}A_2$ defined in [DLMM]. We shall
compose $\tau$ with another automorphism of $V_L$ so that the resulting 
conformal vectors are those in [DLMM] and we can apply 
the decomposition result given in \cite{KMY} for $V_{\sqrt{2}A_2}.$

Let $\varphi$ be the automorphism of $\VL$ defined by 
\[
\varphi : u \ot e^{\be} \longmapsto 
(-1)^{\la -\al_{2}+\al_{3}, \, \be\ra/2} u \ot e^{\be}
\]
for $u \in M(1)$ and $\be \in L$ and set
\[
\rho = (\theta_{2}\ot 1 \ot 1)\varphi\tau,
\]
where $\theta_{2}\ot 1 \ot 1$ is the automorphism which acts as 
$\theta_{2}$ on $V_{\Z\al_{1}}$ and acts as the identity on 
$V_{\Z\al_{2}} \ot V_{\Z\al_{3}}$.  Let $\widet{\om}^{i} = \rho(\om^{i})$. 
\begin{lem} $(1)$ We have 
\begin{align*}
\rho(s^{1}) &= \frac{1}{8}w_{\be_{1}}^{-}, & 
\rho(s^{2}) &= \frac{1}{10}\sum_{\be \in \Phi_{2}^{+}} 
w_{\be}^{-},\\
\rho(s^{3}) &= \frac{1}{6}(\be_{1}(-1)^{2} + \be_{2}(-1)^{2} 
+ (\be_{1}+\be_{2})(-1)^{2}), & 
\rho(\om) &= \om.
\end{align*}

$(2)$ The $\widet{\om}^{1}$, $\widet{\om}^{2}$, and 
$\widet{\om}^{3}$ are the mutually orthogonal conformal 
vectors of $V_{E} \cong V_{\sqrt{2}A_{2}}$ defined in 
\cite{DLMN} and 
\[
\widet{\om}^{4} = \rho(\om) - \rho(s^{3}) = \frac{1}{12}
\gm(-1)^{2}
\]
is the Virasoro element of $V_{F}$ with central charge $1$.
\end{lem}

\prf (2) follows from (1) and the definition of the conformal
vectors in $V_{\sqrt{2}A_2}$ given in (\ref{ea3}). 

(1) follows from Lemma \ref{al1}
and the definitions of all automorphisms involved. For example,
 \begin{align*}
\rho(s^{1}) &= \frac{1}{8} ((\theta_{2}\ot 1 \ot 1)\varphi)w_{\be_{3}}^{+}\\
&=  \frac{1}{16} ((\theta_{2}\ot 1 \ot 1)\varphi)((\alpha_1(-1)-\alpha_2(-1))^2
+4(e^{\alpha_1}e^{-\alpha_2}+e^{-\alpha_1}e^{\alpha_2}))\\
&=  \frac{1}{16} (\theta_{2}\ot 1 \ot 1)((\alpha_1(-1)-\alpha_2(-1))^2
-4(e^{\alpha_1}e^{-\alpha_2}+e^{-\alpha_1}e^{\alpha_2}))\\
&=\frac{1}{16} ((\alpha_1(-1)+\alpha_2(-1))^2
-4(e^{-\alpha_1}e^{-\alpha_2}+e^{\alpha_1}e^{\alpha_2}))\\
&=\frac{1}{8}w_{\be_{1}}^-.
\end{align*}
\qed

Let $\widet{T}'$ be the subalgebra generated by 
$\widet{\om}^{1}$, $\widet{\om}^{2}$, and 
$\widet{\om}^{3}$, and $\widet{T}''$ the subalgebra generated by 
$\widet{\om}^{4}$. Then $\widet{T}' \subset V_{E}^{+}$ and 
$\widet{T}'' \subset V_{F}^{+}$.  Moreover,
\[
  \begin{array}{ll}
\widet{T}' \cong L(\frac{1}{2},\,0) \ot L(\frac{7}{10},\,0) \ot 
L(\frac{4}{5},\,0), \qquad &
\widet{T}'' \cong L(1,\,0).
  \end{array}
\]

As a $\widet{T}'$-module $V_{E}^{\pm}$ decomposes into a 
direct sum of irreducible $\widet{T}'$-submodules of the form 
$L(\frac{1}{2},\,h_{1}) \ot L(\frac{7}{10},\,h_{2}) \ot 
L(\frac{4}{5},\,h_{3})$.  It follows from 
\cite[Lemma 4.1]{KMY} that $V_{E}^{+}$ is a direct sum of 
four irreducible submodules, which are isomorphic to 
\begin{equation}
  \begin{array}{ll}
  \medskip
L(\frac{1}{2},\,0) \otimes L(\frac{7}{10},\,0) \otimes 
L(\frac{4}{5},\,0), \qquad &
L(\frac{1}{2},\,0) \otimes L(\frac{7}{10},\,\frac{3}{5}) \otimes 
L(\frac{4}{5},\,\frac{7}{5}),\\
\medskip
L(\frac{1}{2},\,\frac{1}{2}) \otimes 
L(\frac{7}{10},\,\frac{1}{10}) \otimes L(\frac{4}{5},\,\frac{7}{5}), 
& 
L(\frac{1}{2},\,\frac{1}{2}) \otimes 
L(\frac{7}{10},\,\frac{3}{2}) \otimes L(\frac{4}{5},\,0),
  \end{array}
\label{eq:3.1}
\end{equation}
and $V_{E}^{-}$ is a direct sum of four irreducible submodules, 
which are isomorphic to
\begin{equation}
  \begin{array}{ll}
  \medskip
L(\frac{1}{2},\,0) \otimes L(\frac{7}{10},\,\frac{3}{5}) \otimes 
L(\frac{4}{5},\,\frac{2}{5}), \qquad 
&
L(\frac{1}{2},\,\frac{1}{2}) \ot L(\frac{7}{10},\,\frac{1}{10}) 
\ot L(\frac{4}{5},\,\frac{2}{5}),\\
\medskip
L(\frac{1}{2},\,0) \otimes L(\frac{7}{10},\,0) \otimes 
L(\frac{4}{5},\,3), 
& 
L(\frac{1}{2},\,\frac{1}{2}) \otimes L(\frac{7}{10},\,\frac{3}{2}) 
\otimes 
L(\frac{4}{5},\,3).
  \end{array}
\label{eq:3.2}
\end{equation}

In \cite[Lemma 4.2]{KMY} 
the irreducible $\widet{T}'$-submodules 
with minimal weight $2/3$ of  
$V_{E+\sqrt{2}(\be_{1}-\be_{2})/2}$, which is denoted by  
$V^{2}$ in \cite{KMY},  
are determined. By using fusion 
rules (\cite{DMZ}, \cite{W}) 
we see that as a $\widet{T}'$-module 
$V_{E+\sqrt{2}(\be_{1}-\be_{2})/2}$ is a direct sum of 
four irreducible submodules, which are isomorphic to 
\begin{equation}
  \begin{array}{ll}
\medskip
L(\frac{1}{2},\,0) \otimes L(\frac{7}{10},\,0) \otimes 
L(\frac{4}{5},\,\frac{2}{3}), \qquad & 
L(\frac{1}{2},\,0) \otimes L(\frac{7}{10},\,\frac{3}{5}) \otimes 
L(\frac{4}{5},\,\frac{1}{15}),\\
\medskip
L(\frac{1}{2},\,\frac{1}{2}) \otimes 
L(\frac{7}{10},\,\frac{1}{10}) \otimes L(\frac{4}{5},\,\frac{1}{15}),
& 
L(\frac{1}{2},\,\frac{1}{2}) \otimes 
L(\frac{7}{10},\,\frac{3}{2}) \otimes L(\frac{4}{5},\,\frac{2}{3}).
  \end{array}
\label{eq:3.3}
\end{equation}

The decompositions of $V_{F}^{\pm}$ and $V_{F+\gm/3}$ as 
$\widet{T}''$-modules can be found in \cite{DG}; that is,
\begin{equation}
\begin{split}
V_{F}^{+} & \cong (\oplus_{m \ge 0} L(1,\,4m^{2})) 
\oplus (\oplus_{m \ge 1} L(1,\,3m^{2})),\\
V_{F}^{-} & \cong (\oplus_{m \ge 0} L(1,\,(2m+1)^{2})) 
\oplus (\oplus_{m \ge 1} L(1,\,3m^{2})),\\
V_{F+\gm/3} & \cong \oplus_{m \in \Z} L(1,\,(3m+1)^{2}/3).
\end{split}
\label{eq:3.4}
\end{equation}

{}From these decompositions and Lemma \ref{l3.3} we know all 
irreducible direct summands of $\VL^{+}$ as a 
$\widet{T}' \ot \widet{T}''$-module.

Finally, note that the automorphisms 
$(\theta_{2}\ot 1\ot 1)\varphi$ and $\psi_{2}$ of $\VL$ 
commute. Thus $\rho\psi_{1} = \psi_{2}\rho$ and $\rho(V_{N}) 
= \VL^{+}$.  Since $\rho(T) = \widet{T}' \ot \widet{T}''$, 
the decomposition of 
$V_{N}$ as a $T$-module and 
the decomposition of 
$\VL^{+}$ as a $\widet{T}' \ot \widet{T}''$-module 
are the same. Using \eqref{eq:3.1}, 
\eqref{eq:3.2}, \eqref{eq:3.3}, \eqref{eq:3.4},  
and Lemma \ref{l3.3} we conclude:  

\begin{thm} The decomposition of $V_{N}$ into a direct sum of 
irreducible $T$-submodules is as follows:
\begin{eqnarray*}
& &V_{\sqrt{2}A_3}=\left(L(\frac{1}{2},\,0) \otimes L(\frac{7}{10},\,0) \otimes 
L(\frac{4}{5},\,0)\bigoplus L(\frac{1}{2},\,0) \otimes L(\frac{7}{10},\,\frac{3}{5}) \otimes 
L(\frac{4}{5},\,\frac{7}{5})\right.\\
& &\left.\ \ \ \ \ \ \  \bigoplus L(\frac{1}{2},\,\frac{1}{2}) \otimes 
L(\frac{7}{10},\,\frac{1}{10}) \otimes L(\frac{4}{5},\,\frac{7}{5})\bigoplus 
L(\frac{1}{2},\,\frac{1}{2}) \otimes 
L(\frac{7}{10},\,\frac{3}{2}) \otimes L(\frac{4}{5},\,0)\right)\\
& &\ \ \ \ \ \ \   \ \ \bigotimes \left((\oplus_{m \ge 0} L(1,\,4m^{2})) 
\bigoplus (\oplus_{m \ge 1} L(1,\,3m^{2}))\right)\\
& &\ \ \ \ \ \ \  \bigoplus \left(L(\frac{1}{2},\,0) \otimes L(\frac{7}{10},\,\frac{3}{5}) \otimes 
L(\frac{4}{5},\,\frac{2}{5})\bigoplus
L(\frac{1}{2},\,\frac{1}{2}) \ot L(\frac{7}{10},\,\frac{1}{10}) 
\ot L(\frac{4}{5},\,\frac{2}{5})\right.\\
& &\left.\ \ \ \ \ \ \  \bigoplus 
L(\frac{1}{2},\,0) \otimes L(\frac{7}{10},\,0) \otimes 
L(\frac{4}{5},\,3)\bigoplus 
L(\frac{1}{2},\,\frac{1}{2}) \otimes L(\frac{7}{10},\,\frac{3}{2}) 
\otimes 
L(\frac{4}{5},\,3) \right)\\
& & \ \ \ \ \ \ \   \ \ \bigotimes \left((\oplus_{m \ge 0} L(1,\,(2m+1)^{2})) 
\bigoplus (\oplus_{m \ge 1} L(1,\,3m^{2}))\right)\\
& &\ \ \ \ \ \ \  \bigoplus \left( 
L(\frac{1}{2},\,0) \otimes L(\frac{7}{10},\,0) \otimes 
L(\frac{4}{5},\,\frac{2}{3})\bigoplus 
L(\frac{1}{2},\,0) \otimes L(\frac{7}{10},\,\frac{3}{5}) \otimes 
L(\frac{4}{5},\,\frac{1}{15})\right.\\
& & \left.\ \ \ \ \ \ \  \bigoplus 
L(\frac{1}{2},\,\frac{1}{2}) \otimes 
L(\frac{7}{10},\,\frac{1}{10}) \otimes L(\frac{4}{5},\,\frac{1}{15})
 \bigoplus 
L(\frac{1}{2},\,\frac{1}{2}) \otimes 
L(\frac{7}{10},\,\frac{3}{2}) \otimes L(\frac{4}{5},\,\frac{2}{3})\right)\\
& & \ \ \ \ \ \ \ \ \    \bigotimes \left(\oplus_{m \in \Z} L(1,\,(3m+1)^{2}/3)\right).
\end{eqnarray*}
\end{thm}

\end{document}